\documentclass[12pt]{amsart}

\newtheorem{theorem}{Theorem}[section]
\newtheorem{proposition}[theorem]{Proposition}
\newtheorem{corollary}[theorem]{Corollary}
\newtheorem{remark}[theorem]{Remark}
\newtheorem{lemma}[theorem]{Lemma}
\newtheorem{definition}[theorem]{Definition}
\usepackage{amssymb}

\numberwithin{equation}{section}

\begin{document}

\baselineskip=16pt

\title[Non-emptiness]{Non-emptiness of moduli
 spaces of coherent systems}

\author{L. Brambila-Paz}

\address{CIMAT, Apdo. Postal 402, C.P. 36240.
 Guanajuato, Gto,
M\'exico}

\email{lebp@cimat.mx}

\keywords{coherent systems, moduli space, stability,
Brill-Noether}

\subjclass[2000]{14H60, 14J60}

\date{}
\begin{abstract}
Let $X$ be a general smooth projective algebraic curve of genus
$g\geq 2$ over $\mathbb{C}$. We prove that the moduli space
$G(\alpha:n,d,k)$ of $\alpha $-stable coherent systems of type
$(n,d,k)$ over $X$ is empty if $k>n$ and the Brill-Noether number
$\beta:=\beta(n,d,n+1)=\beta(1,d,n+1)= g-(n+1)(n-d+g)<0$.
Moreover, if $0\leq \beta <g$ or $\beta =g, n\not|g$ and for some
$\alpha
>0$, $G(\alpha:n,d,k)\not= \emptyset$ then $G(\alpha :n,d,k)\not=
\emptyset$ for all $\alpha
>0$ and $G(\alpha:n,d,k)= G(\alpha ':n,d,k)$ for all $\alpha
,\alpha '>0$ and the generic element is generated. In particular,
 $G(\alpha:n,d,n+1)\not= \emptyset$ if $0\leq \beta \leq g$ and $\alpha>0$.
 Moreover, if $\beta >0$ $G(\alpha:n,d,n+1)$ is
 smooth and irreducible of dimension $\beta(1,d,n+1).$ We define a
 dual span of a generically generated coherent
system. We assume $d< g+n_1\leq g+n_2$ and prove that for all
$\alpha
 >0$, $G(\alpha:n_1,d, n_1+n_2)\not= \emptyset $ if and only if
$G(\alpha:n_2,d, n_1+n_2)\not= \emptyset .$ For $g=2$, we describe
 $G(\alpha:2,d,k)$ for $k>n$.

\end{abstract}

\maketitle

\section{Introduction}

Let $X$ be a  smooth projective algebraic curve of genus $g\geq 2$
over $\mathbb{C}$.
 A coherent system over $X$ of
type $(n,d,k)$  is a pair $(E,V)$ where $E$ is a vector bundle
over $X$ of rank $n$, degree $d$ and $V$ a linear subspace of
$H^0(X,E)$ of dimension $k$.

A notion of stability for coherent systems was introduced in
\cite{lep}, \cite{rag} and \cite{an}. The definition of stability
depends on a real parameter $\alpha $, which corresponds to the
choice of linearization of a group action. The
  coherent systems are also ``augmented bundles''
(see \cite{bdo}) and are related with the existence of solutions
of orthogonal vortex equations, where the parameter $\alpha$
appears in a natural way.

For any $\alpha \in \mathbb{R}$ denote by $G(\alpha:n,d,k)$
(respectively $\widetilde{G}(\alpha:n,d,k)$) the moduli space of
$\alpha $-stable (respectively $\alpha$-semistable) coherent
systems of type $(n,d,k).$ From the definition of
$\alpha$-stability one can see that in order to have
$\alpha$-stable coherent systems with $k\geq 1$ we need $\alpha
>0$. The expected dimension of $G(\alpha:n,d,k)$ is the Brill-Noether
number $\beta (n,d,k):=n^2(g-1)+1-k(k-d+n(g-1)).$  Note that if
$k>n$, $\beta (n,d,k)=\beta(k-n,d,k).$ We denote by $\beta$ the
Brill-Noether number $\beta(n,d,n+1)=\beta(1,d,n+1)=
g-(n+1)(n-d+g)$.

 Basic properties
of $G(\alpha:n,d,k)$ have been proved in \cite{lep}, \cite{an},
\cite{rag} and particular cases in \cite{bu}, \cite{bo} and
\cite{bommn}. More general results can be found in \cite{bomn},
\cite{he} and \cite{bdo}. Most of the detailed results known are
for $k\leq n.$ It is our purpose here to study the case $k>n.$

In \cite[Proposition 4.6]{bomn} it was proved that for $k\geq n$
there exists $\alpha _L$ such that $G(\alpha:n,d,k)= G(\alpha '
:n,d,k)$ if $\alpha , \alpha ' >\alpha _L.$ Denote this moduli
space by $G_{L}(n,d,k).$

 For any $(n,d,k)$ define $U(n,d,k)$ and
$U^s(n,d,k)$ as
$$U(n,d,k):=\{(E,V):(E,V)\in G_L(n,d,k)\mbox{ and }E\mbox{ is stable}\}$$
and
$$
U^s(n,d,k):= \{ (E,V): (E,V)\  {\rm is \  of \  type }\  (n,d,k) \
 {\rm and \   is} \  \alpha {\rm -stable \  for \  all} \ \
\alpha >0 \}.
$$

We prove the following (see Theorem \ref{t1})

\noindent \textbf{Theorem 1}\, \textit{ Let $X$ be general, $\beta
<g$ or $\beta =g, n\not|g$ and $k> n$. Then
\begin{enumerate}
\item if $\beta <0$, $G(\alpha:n,d,k)=\emptyset $ for all $\alpha
>0$; \item if for some $\alpha >0$, $G(\alpha:n,d,k)\not=
\emptyset$ then $G(\alpha:n,d,k)\not= \emptyset$ for all $\alpha
>0$ \item \textit{$G(\alpha:n,d,k)= G(\alpha ':n,d,k)$ for all
$\alpha ,\alpha '>0$ i.e. $\alpha _L=0$};
 \item $(E,V) \in G(\alpha: n,d,k)$ if and only if $(E,V)$
 is generically generated and $H^0(I_E^*)=0$, where
 $I_E$ is the image of the evaluation map $V\otimes
\mathcal{O}\rightarrow E ;$
\item if for some $\alpha >0$,
$G(\alpha:n,d,k)\not= \emptyset$ then $U(n,d,k)= G(\alpha:n,d,k)$.
\end{enumerate} }

Note that the results of {\it Theorem 1} deal with the moduli
spaces of coherent systems of type $(n,d,k)$ whereas $\beta$
refers to $(n,d,n+1)$. Moreover, if $\beta (n,d,n+1) \leq g$,
$\beta(n,d,k) <0$ for $k>n+1.$

 If $\alpha _L =0$ denote $G_L(n,d,k)$
by $G(n,d,k)$. In particular,  $\mathcal{G}^{k-1}_d:=G(1,d,k).$
For $k=n+1$ we have (see Theorem \ref{t21})

\noindent \textbf{Theorem 2}\, \textit{ Let $X$ be general and $
\beta := \beta(n,d,n+1)\leq g$. Then
\begin{enumerate}
\item $G(\alpha:n,d,n+1)\not= \emptyset$ if and only if $\beta
\geq 0$; \item if $\beta \geq 0$ then $G(n,d,n+1):= G(\alpha
:n,d,n+1) = G(\alpha ':n,d,n+1)$ for all $\alpha ,\alpha ' >0$ and
$\alpha _L =0;$ \item if $\beta> 0$ then $G(n,d,n+1)$ is smooth
and irreducible of dimension $ \beta $ and the generic element is
generated; \item $U^s(n,d,n+1)= G(n,d,n+1)$ and is birationally
equivalent to $\mathcal{G}^{n}_d$; \item if $\beta =0$
$G(n,d,n+1)\cong \mathcal{G}^{n}_d$ and the number of points of
$G(n,d,n+1)$ is
$$g!\prod_{i=0}^n\frac{i!}{(g-d+n+i)!}.$$
\end{enumerate} }

Moreover, (see Theorem \ref{t23})

\noindent \textbf{Theorem 3}\, \textit{ If $X$ is general and
$g\geq n^2-1$ then for any degree $d\geq g+n-\frac{g}{n+1}$
\begin{enumerate}
\item $G(\alpha:n,d,n+1)\not= \emptyset$ for all $\alpha > 0$;
\item $U(n,d,n+1)\not= \emptyset $ and is smooth and irreducible.
\end{enumerate}}

 As was pointed out in \cite{bo} and \cite{bomn}
coherent systems are related with Brill-Noether theory.  Let
$B(n,d,k)$ (respectively $\widetilde{B}(n,d,k)) $ be the
Brill-Noether locus defined by stable (respectively semistable)
 vector bundles of rank $n$,
degree $d$ and $\dim H^0(X,E)\geq k$. It is well known that for
``small" $\alpha $, $(E,V)$ $\alpha$-stable implies $E$ semistable
and $E$ stable implies $(E,V)$ $\alpha $-stable. The approach to
study the Brill-Noether loci in \cite{bomn} is to describe
$G(\alpha:n,d,k)$, usually for ``large" $\alpha $, and through
``flips" obtain information of $G(\alpha:n,d,k)$ for smaller
$\alpha . $

In our case, i.e. $\beta <g$ or $\beta =g, n\not|g$ and $k> n$, it
is enough to know the non-emptiness for one $\alpha$ to obtain
non-emptiness for all $\alpha.$ Moreover, there are no ``flips".

 In \cite{ye} it was proved
that if $X$ is general and $g\geq \beta (n,d,n+1) \geq 0 $,
$B(n,d,n+1)$ is non-empty and has a component of the correct
dimension. From the above results of coherent systems we have (see
Corollary \ref{c21})

\noindent \textbf{Corollary 4}\, \textit{ If $X$ is general and
$g\geq \beta\geq 0$, $B(n,d,n+1)$ is irreducible if $\beta
>0$ and $G(n,d,n+1)$ is a desingularisation of (the closure of) the
Brill-Noether locus $B(n,d,n+1).$ Moreover, the natural map $\phi
: G(\alpha:n,d,n+1) \rightarrow \widetilde{B}(n,d,n+1)$ is an
isomorphism  on the complement of the singular locus of
 $B(n,d,n+1)\subset \widetilde{B}(n,d,n+1)$.}

Actually, \cite[Conditions 11.3]{bomn} are satisfied in this case
and hence the results in \cite[\S 11]{bomn} hold.

Besides the known relation between coherent systems and
Brill-Noether theory, our results on $G(n,d,n+1)$ can be related
with other problems. Given a generated linear system $(L,V)$ we
have the natural map $$\phi _V :X\rightarrow \mathbb{P}(V^*).$$ In
particular, if $L$ has degree $d$ and dim$V=n+1$ we have (see
Theorem \ref{t3});

\noindent \textbf{Theorem 5}\, \textit{ Let $X$ be general, $0
\leq \beta (n,d,n+1)$ and $T\mathbb{P}$ the tangent bundle of
$\mathbb{P}(V^*).$ If $\beta<g$ or $\beta=g$ and $n\not|g$, then
$\phi_V^*(T\mathbb{P})$ is stable. If either $g\ge n^2-1$ or
$\beta=g$, $n|g$ and $g$ and $n$ are not both equal to $2$, then
there exist linear systems $(L,V)$ such that
$\phi_V^*(T\mathbb{P})$ is stable.}

We define a dual span of a generically generated coherent system
(see Definition \ref{def1}) and denote by $D(E,V)=
(D(E)_\ell,V^*)$ a dual span of $(E,V)$. If $I_E$ is the image of
the evaluation map $V\otimes \mathcal{O}\rightarrow E$ we prove
(see Theorems \ref{t31} and \ref{theorem1})

\noindent \textbf{Theorem 6}\, \textit{Let $X$ be a general
 curve of genus $g$ and $d< g+n_1\leq g+n_2$ then for all $\alpha
 >0$, $G(\alpha:n_1,d, n_1+n_2)\not= \emptyset $ if and only if
$G(\alpha:n_2,d, n_1+n_2)\not= \emptyset $.}

 \noindent\textbf{Theorem 7}\, \textit{ Let $(E,V) \in
G(\alpha:n_1,d,n_1+n_2)$. If either of the Petri maps of $(I_E,V)$
or $(I_{{D(E)_\ell}},V^*)$ is injective then,
\begin{enumerate} \item $G(\alpha:n_1,d,n_1+n_2)$ is smooth of dimension
$\beta(n_1,d,n_1+n_2)$ in a neighbourhood of $(E,V)$; \item
$G(\alpha:n_2,d,n_1+n_2)$ is smooth of dimension
$\beta(n_2,d,n_1+n_2)$ in a neighbourhood of the dual span
$D(E,V).$
\end{enumerate}}

Denote by $G_0(n,d,k)$ the moduli space $G(\alpha:n,d,k)$ for
''small" values of $\alpha $ (see Remark \ref{rem1} (2)). For
$n=2$ we have (see Theorem \ref{teo7})

 \noindent \textbf{Theorem 8}\, \textit{Let $X$ be general, $s\geq 3$
and $d< s +2g -\frac{4g}{s+2}$. If $G_0(2,d,2+s)$ is non-empty
then $G(\alpha :2,d,2+s)$ is non-empty for all $\alpha >0.$
Moreover, $U(2,d,2+s)\not= \emptyset $and $g=2 $}

For $n=2$ and $g=2$, from the above results and the Riemann-Roch
Theorem we know that

\begin{enumerate} \item if $d<4$ and $k\geq
3$, $G(\alpha :2,d,k)=\emptyset$ for all $\alpha
>0$;
 \item if $d= 5$ and $k> 3$, $U(2,d,k)=\emptyset $ and $G_0(2,5,k)=\emptyset$; \item if
$d\geq 6$ and $k=3,4$,  $G(\alpha :2,d,k)\not =\emptyset$ for all
$\alpha >0$. Moreover, $U(2,d,k)\not=\emptyset ;$ \item if $d \geq
6$ and $k> d-2$, $U(2,d,k)=\emptyset $ and $G_0(2,d,k)=\emptyset
.$
\end{enumerate}

In particular for $d=4,5$ we have (see Theorems
\ref{proposition4}, \ref{teoff} and \ref{prop5})

\noindent \textbf{Theorem 9}\, \textit{ \begin{enumerate}\item
$U(2,4,k)=\emptyset $ for $k\geq 3$; \item $G_0(2,4,k)=\emptyset$
for $k\geq 5;$ \item $G(\alpha :2,4,3)\not= \emptyset $ for all
$\alpha >0$; \item $U^s(2,4,3)\cong G_L(2,4,3)\cong Pic^4(X)$.
\end{enumerate}}

\noindent \textbf{Theorem 10}\, \textit{\begin{enumerate} \item
$\widetilde{G}_0(2,4,4)=\{(K\oplus K , H^0(K\oplus K))\}$; \item
$G_0(2,4,4)= \emptyset$;  \item $\widetilde{G}(\alpha:2,4,4)\not=
\emptyset$ for all $\alpha >0$.
\end{enumerate}}

\noindent \textbf{Theorem 11}\, \textit{
\begin{enumerate} \item $G_0(2,5,3)\not= \emptyset$;
 \item $U(2,5,3)\not =\emptyset $; \item $U(2,5,3)\not=
 G_0(2,5,3).$
\end{enumerate} }

{\it Notation}

We will denote by $K$ the canonical bundle over $X$, by $I_E$ the
image of the evaluation map $V\otimes \mathcal{O}\rightarrow E$,
 $H^i(X,E)$ by $H^i(E)$, $\dim H^i(X,E)$
by $h^i(E)$, the rank of $E$ by $n_E$, the degree of $E$ by $d_E$
and $\det (E)$ by $L_E$. By a general curve we mean a Petri curve
i.e. the Petri map
$$H^0(L)\otimes H^0(L^*\otimes K)\rightarrow H^0(K)$$
is injective for every  line bundle $L$ over $X$.

\bigskip

\section{General results}

Let $X$ be an irreducible smooth projective curve over $\mathbb C$
of genus $g\geq 2$.
 For any
 $\alpha \in \mathbb{R},$
define the $\alpha$-slope of the coherent system $(E,V)$ of type
$(n,d,k)$ as
$$\mu _{\alpha}(E,V):= \mu(E) + \alpha \frac{k}{n},$$
where $\mu (E):=d/n$ is the slope of the vector bundle $E.$ A
coherent subsystem $(F,W)\subseteq (E,V)$ is a
 coherent system such that $F\subseteq E $ and $ W\subseteq V\cap H^0(F).$
For any $\alpha \in \mathbb{R}$ a coherent system $(E,V)$ is
$\alpha$-stable (respectively $\alpha$-semistable) if for all
proper coherent subsystems $(F,W)$
$$\mu _{\alpha}(F,W)<\mu _{\alpha}(E,V) \ \ \
({\rm respectively} \leq ).$$

Denote the moduli space of $\alpha$-stable (respectively $
\alpha$-semistable) coherent systems of type $(n,d,k)$ by
$G(\alpha:n,d,k)$ (respectively $\tilde{G}(\alpha:n,d,k)$) and by
$\beta (n,d,k) $ the Brill-Noether number $\beta (n,d,k):=
n^2(g-1)+1 -k(k-d+n(g-1)).$ From the infinitesimal study of
 the coherent systems (see
\cite{bomn} and \cite{he}) we have that

\begin{proposition}\label{e} If $(E,V)\in G(\alpha:n,d,k)$
then $G(\alpha:n,d,k)$ is smooth of dimension $\beta (n,d,k) $
 in a neighbourhood of $(E,V)$ if
and only if the Petri map $V\otimes H^0(E^*\otimes K)\rightarrow
H^0(E\otimes E^*\otimes K)$ is injective.
 Moreover, $ T_{(E,V)} G(\alpha:n,d,k)=
Ext^1((E,V),(E,V)).$
\end{proposition}

If $B(n,d,k)$ (respectively $\widetilde{B}(n,d,k)$) is the
Brill-Noether locus of stable (respectively semistable) vector
bundles then for ``small'' $\alpha $ there is a natural map
$$\phi :  G(\alpha:n,d,k) \rightarrow  \widetilde{B}(n,d,k)$$
defined by $(E,V) \mapsto E$ that is injective over
$B(n,d,k)-B(n,d,k+1).$

Given a triple $(n,d,k)$ denote by $C(n,d,k)$ the set
$$C(n,d,k):=\{\alpha \in \mathbb{R}| 0\leq \alpha
=\frac{nd'-n'd}{n'k-nk'} \  {\rm with} \ \ 0\leq k'\leq k,
0<n'\leq n, \ {\rm and} \ \
 nk'\not= n'k \}.$$
An element $\alpha $ in $ C(n,d,k)$ is
 called a virtual
critical point. The set $C(n,d,k)$  defines a
 partition of the interval $[0,\infty )$. With
 the natural
order on $\mathbb{R}$, label the virtual critical points as
$\alpha _i.$

It is known (see \cite{bdo} and \cite{bomn}) that

\begin{remark}\begin{em}\label{rem1}
\begin{enumerate}
\item If $(n,d,k)=1$ then $G(\alpha:n,d,k)=
\widetilde{G}(\alpha:n,d,k)$, for $\alpha \not\in C(n,d,k).$ \item
If $\alpha ',\alpha '' \in (\alpha _i , \alpha _{i+1})$ then
  $G(\alpha ':n,d,k)= G(\alpha '':n,d,k)$.
Denote by
 $G_i(n,d,k)$ the moduli space $G(\alpha:n,d,k)$
for any $\alpha \in (\alpha _i , \alpha _{i+1}).$ \item For $k\geq
n,$ there exists $\alpha _L$ such that for any $\alpha , \alpha
'>\alpha _L, \ \ G(\alpha:n,d,k)= G(\alpha ':n,d,k)$. Denote by
$G_L(n,d,k)$ the moduli space
 $G(\alpha:n,d,k)$ for $\alpha > \alpha _L$.
\item Every irreducible component of $G_i(n,d,k)$
 has dimension at least
$\beta (n,d,k).$
\end{enumerate}
\end{em}
\end{remark}

\begin{remark}\begin{em}\label{r0} Let $(E,V)$ be a coherent system of type
$(n,d,k)$. From the definition of $\alpha $-stability and
stability of a vector bundle we have that
\begin{enumerate}
\item if $(E,V)\in G(\alpha:n,d,k)$ and $E$ is stable then $(E,V)$
is $\alpha '$-stable for all $0<\alpha ' <\alpha ;$ \item if $E$
is stable and for all coherent subsystems $(F,W) \subset (E,V)$,
$\frac{{\rm dim}W}{n_F} \leq \frac{k}{n}$ then $(E,V)$ is
$\alpha$-stable for all $\alpha >0;$ \item if $E$ is semistable
and for all coherent subsystems $(F,W) \subset (E,V)$, $\frac{{\rm
dim}W}{n_F} < \frac{k}{n}$ then $(E,V)$ is $\alpha$-stable for all
$\alpha >0;$ \item if $E$ is semistable and for all coherent
subsystems $(F,W)\subset(E,V)$, $\frac{\dim
W}{n_F}\le\frac{k}{n}$, then $(E,V)$ is $\alpha$-semistable for
all $\alpha >0 .$

\end{enumerate}\end{em}
\end{remark}

Let $(E,V)$ be a coherent system of type $(n, d,k)$ with $k> n.$
We shall say that $(E,V)$ (or $E$) is generically generated if the
image $I_E$ of the evaluation map $V\otimes \mathcal{O}\rightarrow
E$ has rank $n.$ That is, we have the exact sequence
\begin{equation}\label{gge}
0\rightarrow I_E \rightarrow E \rightarrow \tau \rightarrow 0
\end{equation}
where $\tau $ is a torsion sheaf.  If $\tau =0$
 we say
that $(E,V)$ (or $E$) is generated and if $\tau \not= 0$, is
strictly generically generated.

\begin{remark}\begin{em}\label{pp11} Note that if $(E,V)$ is generated with
$H^0(E^*)=0$, any quotient bundle $Q$ is generated and
$H^0(Q^*)=0.$
\end{em}\end{remark}

We give a proposition that we will use in the following sections

\begin{proposition}\label{p21} Let $(E,V)$ be a generated coherent system of
type $(n,d,k)$ with $E$ semistable and $k=n+s, s\geq 1$. If
$(F,W)$ is a coherent subsystem of $(E,V)$:
\begin{enumerate}
\item  $\dim W \leq n_F+s-1;$ \item  if $\frac{(s-1)n}{s}< n_F$,
$\mu _{\alpha} (F,W) < \mu_{\alpha} (E,V)$ for all $\alpha
>0$; \item if $\dim W\leq n_F$, $\mu _{\alpha} (F,W)< \mu_{\alpha} (E,V)$
for all $\alpha>0$;  \item if $(E,V)$ is of type $(n,d,n+1) $ then
it is $\alpha $-stable for all $\alpha >0.$
\end{enumerate}
\end{proposition}

\begin{proof} Note that $d>0$, so $E^*$ is semistable of
negative degree, hence $H^0(E^*)=0.$ Let $(F,W)$ be a coherent
subsystem of $(E,V)$ and $(Q,Z)$ the quotient coherent system.
Since $Q$ is generated and $H^0(Q^*)=0$,
$$ \dim (V)-\dim (H^0(F)\cap V) \geq n_Q+1$$
that is, $n_F +s-1\geq dim(H^0(F)\cap V)\geq \dim W$.

If $\frac{(s-1)n}{s}< n_F$, $\frac{\dim(W)}{n_F}<
\frac{\dim(V)}{n}$ and from Remark $2.3$
$$\mu _{\alpha} (F,W) < \mu_{\alpha} (E,V)$$
for all $\alpha >0.$ Similarly, for $\dim W\leq n_F$, $\mu
_{\alpha} (F,W) < \mu_{\alpha} (E,V)$ for all $\alpha >0$.

If $s=1$, for all coherent subsystems $(F,W)$, $n_F \geq \dim W$,
therefore, from Remark $2.3$, $(E,V)$ is $\alpha $-stable for all
$\alpha
>0.$
\end{proof}

For any $(n,d,k)$ define $U^s(n,d,k)$ and $U(n,d,k)$ as

\begin{equation}\label{e3}
U^s(n,d,k):= \{ (E,V): (E,V)\  {\rm is \  of \  type }\  (n,d,k) \
 {\rm and \   is} \  \alpha {\rm -stable \  for \  all} \ \
\alpha >0 \};
\end{equation}

$$ U(n,d,k):=\{(E,V):(E,V)\in G_L(n,d,k)\mbox{ and }E\mbox{ is
stable}\}$$

From Remark \ref{r0}(1) we have that $U(n,d,k)\subset U^s(n,d,k).$
Note that $U^s(n,d,k)$ is embedded in $G_L(n,d,k).$ From the
openness of $\alpha$-stability it follows that $U^s(n,d,k)$ is an
open subset of $G_L(n,d,k)$. Moreover, if $(E,V)\in U^s(n,d,k)$,
$E$ is semistable.

\begin{proposition}\label{p1dg} If $d\geq n(2g-1)$,
$G(\alpha:n,d,n+1)\not= \emptyset$ for all $\alpha >0$. Moreover,
$U(n,d,n+1)\not= \emptyset .$
\end{proposition}

\begin{proof} If $d\ge n(2g-1)$, every stable bundle $E$ of rank $n$
and degree $d$ is generated and $h^0(E)\ge n+1$. A generic
subspace $V$ of $H^0(E)$ of dimension $n+1$ generates $E$. By
Proposition \ref{p21}(4), $(E,V)$ is $\alpha$-stable for all
$\alpha>0$.
 Hence $U(n,d,n+1)\ne\emptyset$.
\end{proof}

Our aim is to prove that such coherent systems exist for smaller
$d$.

\section{Vector bundles with sections}

In this section we assume that $X$ is a general curve and $k\geq
n+1$. We give three Lemmas that we will use.

\begin{lemma}\label{re1} If $F$ is generated and $H^0(F^*)=0,$ then
$\mu (F) \geq 1+\frac{g}{n_F+1}.$
\end{lemma}

\begin{proof}
Recall from \cite[Proposition 3.2]{pr} that if $F$ is generated
and $H^0(F^*)=0$ then it
 is generated by a linear subspace $W\subseteq H^0(F)$ of dimension
$n_F+1$, and $h^0(det (F))\geq n_F+1.$ Moreover, the Brill-Noether
theory for line bundles implies that
$$ \beta(1,d_F,n_F+1)= g-(n_F+1)(n_F-d_F+g)\geq 0.$$
That is,
\begin{equation}\label{b5}
\mu(F)\geq 1+\frac{g}{n_F+1}.
\end{equation}
\end{proof}

\begin{lemma}\label{p1} Let $E$ be a vector
bundle such that $d_E \leq n_E +g.$ If $F$ is a vector bundle of
rank $n_F<n_E$ that is generically generated and $H^0(I_F^*)=0$
then $\mu (F)\geq \mu (E)$. Moreover, $\mu (F)=\mu(E)$ is possible
 only if $n_F=n_E-1.$
\end{lemma}

\begin{proof} By hypothesis $\mu (E) \leq \frac{g}{n_E} +1.$ If $\mu (I_F) \leq \mu
(F)<\mu(E)$ then from Lemma \ref{re1} we get a contradiction.
\end{proof}

\begin{corollary}\label{c1} If $E$ is a semistable bundle with
$d_E< n_E+g$ or $d_E=g+n_E$, $n_E\not|g$ then $E$ can not have a
proper generically generated subbundle $F$ with $H^0(I_F^*)=0$.
\end{corollary}

\begin{proof} Suppose that $F\subset E$ is generically generated with $H^0(I_F^*)=0$. From the
semistability and Lemma $3.2$, $n_F=n_E-1$ and $d_E=g+n_E$. But
then $E/F$ is a line bundle and $\mu (F)=\mu(E)=\mu(E/F)$ which is
a contradiction if $d_E=g+n_E$, $n_E\not|g$.

\end{proof}

\begin{lemma}\label{l1} If $F$ is generated by a subspace $W$ of
dimension dim$W\geq n_F+1$ then either $H^0(F^*)=0$ or there is a
subbundle $G$ with $n_G<n_F$ that is generated and $H^0(G^*)=0$.
\end{lemma}

\begin{proof} If $H^0(F^*)\not= 0$ then $F\cong
\mathcal{O}^s \oplus G$ where $G$ is generated, $H^0(G^*)=0$ and
$1\leq n_G<n_F$.
\end{proof}

For coherent systems of type $(n,d,k) $ with $k\geq n+1$ we have
the following propositions.

\begin{proposition}\label{p3} Let $(E,V)$ be a coherent system
of type $(n,d,k)$ with $d < n+g$ or $d=g+n$, $n\not|g.$ Then $E$
is stable if and only if $(E,V)$ is generically generated and
$H^0(I_E^*)=0.$ Moreover, if $d=g+n$, $n|g$ and $(E,V)$ is
generically generated with $H^0(I_E^*)=0$, $E$ is semistable.
\end{proposition}

\begin{proof} Suppose $E$ is stable. Then $I_E$
is generated by $V$. If $H^0(I_E^*)=0$, from Corollary \ref{c1},
$n_{{I_E}}=n_E.$ If $H^0(I_E^*)\not=0$, from Lemma \ref{l1} and
Corollary \ref{c1} we get a contradiction

Now suppose $(E,V)$ is generically generated with $H^0(I^*_E)=0$.
If $E$ is not stable, let $Q$ be a quotient bundle such that $\mu
(Q)\leq \mu (E)$. We have the following diagram

\begin{equation}
\begin{array}{ccccccccc}
0&\rightarrow &I_E &
\rightarrow &E &\rightarrow &\tau &\rightarrow &0\\
&&\downarrow&&\downarrow&&\downarrow&&\\
0&\rightarrow&Q_1&\rightarrow&Q&\rightarrow&
\tau '&\rightarrow&0\\
&&\downarrow&&\downarrow&&&&\\
&&0&&0&&&&\\
\end{array}
\end{equation}
where $Q_1$ is a quotient bundle of $I_E$ such that $\mu (Q_1)
\leq \mu(Q) ,\ \ n_{{Q_1}}= n_Q$ and since $I_E$ is generated and
$H^0(I_E^*)=0$, $Q_1$ is generated and $H^0(Q_1^*)=0$. Thus,

\begin{equation}\label{e1}
1+\frac{g}{n_{{Q}}+1}\leq \mu(Q_1)\leq \ \mu(Q)\leq \mu(E)
=\frac{d}{n} \leq 1+\frac{g}{n}.
\end{equation}

If $n_Q+1<n$ we get a contradiction. If $n_Q+1=n$, $\mu
(Q)=\mu(E)$ and hence $E$ is semistable. But in that case there
exists a line bundle $L_0$ such that $Q\cong E/L_0$ and $\mu
(E)=\mu (Q)=\mu (L_0)$. This will be a contradiction if $n\not|g.$
Therefore $E$ is stable.
\end{proof}

\begin{proposition}\label{p6} A generically generated coherent system $(E,V)$ of type
$(n,d,k)$ with $d< g+n$ or $d=g+n, n\not|g$ and $H^0(I_E^*)=0$ is
$\alpha $-stable for all $\alpha
>0$.
\end{proposition}

\begin{proof} From Proposition \ref{p3}, $E$ is
stable. Let $(F,W) \subset (E,V)$ be a coherent subsystem of
$(E,V)$ with $n_F<n_E$. If dim$(W)\geq n_F+1$, the evaluation map
defines a subbundle $F'$, with $n_{F'}\leq n_F<n_E$ which is
generically generated with $H^0(F'^*)=0$. From Lemmas \ref{l1} and
$3.2$, $\mu (F')\geq \mu(E)$ which contradicts stability of $E$.
Hence, dim$W\leq n_F$ and from Remark $2.3$, $(E,V)$ is
$\alpha$-stable for all $\alpha
>0$.
\end{proof}

For $k=n+1$ we have

\begin{proposition}\label{newcor1}
A generically generated coherent system $(E,V)$ of type
$(n,d,n+1)$ with $d\le g+n$ and $H^0(I_E^*)=0$ is $\alpha$-stable
for all $\alpha>0$.
\end{proposition}

\begin{proof} From Proposition \ref{p3}, $E$ is semistable.
Let $(Q,W)$ be a proper quotient coherent system of $(E,V)$. Then
$(Q,W)$ is generically generated. Moreover, since $I_Q$ is a
quotient of $I_E$, $H^0(I_Q^*)=0$ and hence $\dim W\ge n_Q+1$. So
$\frac{n+1}n<\frac{\dim W}{n_Q}$ and the result follows from
Remark \ref{r0}(3).
\end{proof}

Conversely

\begin{proposition}\label{p7} If $(E,V)$ is an $\alpha$-stable
coherent system of type $(n,d,k)$ with $d\leq g+n$ then $(E,V)$ is
generically generated and $H^0(I_E^*)=0$. Moreover, $E$ is
semistable and stable if $d < n+g$ or $d=g+n$, $n\not|g.$
\end{proposition}

\begin{proof} Suppose that $I_E=\mathcal{O}^s\oplus G$ with $0\leq s\leq n_{I_E}-1$, $G$
generated, $H^0(G^*)=0$ and $\mu (G)\geq \frac{g}{n_G+1} +1$. From
the $\alpha $-stability of $(E,V)$ we have
$$\mu_{\alpha}(G,H^0(G)\cap V) < \mu_{\alpha }(E,V).$$ That is,
$$ \alpha (\frac{k-s}{n_G} -\frac{k}{n}) < \mu
(E)-\mu(G).$$

If $n_G <n$, then $\mu (E)-\mu(G)\leq \frac{g}{n} +1
-(\frac{g}{n_G+1}+1)\leq 0$, hence
$$ \alpha (\frac{k-s}{n_G} -\frac{k}{n}) < 0$$
which is a contradiction since $s\leq n-n_G$. Hence $n_{I_E}=n$,
$(E,V)$ is generically generated and $H^0(I_E^*)=0$. The last part
follows from Proposition \ref{p3}.
\end{proof}

From Propositions \ref{p3}, \ref{p6} and \ref{p7} we have Theorem
1.

\begin{theorem}\label{t1} Let $X$ be general, $\beta =\beta(n,d,n+1) <g$ or
$\beta =g, \ n\not|g$ and $k\geq n+1$. Then
\begin{enumerate}
\item if $\beta <0$, $G(\alpha:n,d,k)=\emptyset$ for all $\alpha
>0$; \item if for some $\alpha >0$, $G(\alpha:n,d,k)\not=
\emptyset$ then $G(\alpha:n,d,k)\not= \emptyset$ for all $\alpha
>0$; \item \textit{$G(\alpha:n,d,k)= G(\alpha ':n,d,k)$ for all
$\alpha ,\alpha '>0$ i.e. $\alpha _L=0$};
 \item $(E,V) \in G(\alpha: n,d,k)$ if and only if $(E,V)$
 is generically generated and $H^0(I_E^*)=0$;
\item if for some $\alpha >0$, $G(\alpha:n,d,k)\not= \emptyset$
then $U^s(n,d,k)= G(\alpha:n,d,k)$ and $U(n,d,k)\not= \emptyset$.
\end{enumerate}
\end{theorem}

\begin{proof} Recall from the definition of $\beta$ that
$\beta(n,d,n+1)=\beta(1,d,n+1)= g-(n+1)(n-d+g).$ Hence,
$$0\leq \beta \Longleftrightarrow \frac{g}{n+1} +1 \leq
\frac{d}{n}.$$ Moreover,
$$\beta \leq g \Longleftrightarrow d\leq g+n.$$
 If
$(E,V) \in G(\alpha:n,d,k)$, $E$ is generically generated and
$H^0(I^*_E)=0$ (see Proposition \ref{p7}). Hence, by Lemma $3.1$,
$\mu (E) \geq \frac{g}{n+1}+1$ i.e $\beta(n,d,n+1)\geq 0$. Parts
2), 3), 4) and 5) follow from Propositions \ref{p6} and \ref{p7}.
\end{proof}

\begin{corollary}\label{corfinal} If $d<g+n$ and $g\leq n$,
$G(\alpha:n,d,k)=\emptyset$ for all $\alpha >0$ and $k\geq n+1$.
\end{corollary}

\begin{proof} It follows from Theorem \ref{t1} since
the Brill-Noether number is negative.
\end{proof}

\section{Coherent systems of type $(n,d,n+1)$}

From Remark \ref{rem1} we have that $G(\alpha:n,d,n+1)=
\widetilde{G}(\alpha:n,d,n+1)$, for $\alpha \not\in C(n,d,n+1).$

For $d\geq n(2g-1),$ from Proposition \ref{p1dg}, $U(n,d,n+1)\not=
\emptyset .$ For small values of $d$ we have the following
Proposition (see also \cite{bu} and \cite{ye})

\begin{proposition}\label{p22} If $X$ is general and $0\le\beta\le g$, then
\begin{enumerate}
\item there exist generated coherent systems $(E,V)$ with $E$
semistable and in particular $U^s(n,d,n+1)\ne\emptyset$; \item
except when $g=n=2$ and $d=4$, there exist generated coherent
systems $(E,V)$ with $E$ stable and in particular
$U(n,d,n+1)\ne\emptyset$.
\end{enumerate}
\end{proposition}

\begin{proof} (1) The dimension of the subvariety consisting of
line bundles $L$ for which $L$ is not generated by a subspace
$V\subset H^0(L)$ of dimension $n+1$ has dimension
$g-(n+1)(n-(d-1)+g)+1<\beta$, since they define a line bundle of
degree $d-1$ with $n+1$ sections. Thus, from the Brill-Noether
theory for line bundles, the set of generated line bundles $L$ of
degree $d$ with $n+1\leq \dim V \leq h^0(L)$ defines a non-empty
open set of the Jacobian $J^d(X).$

 We have the following exact sequence
\begin{equation}\label{e2}
 0\rightarrow E^* \rightarrow
V\otimes \mathcal{O} \rightarrow L \rightarrow 0.
\end{equation}

The coherent system $(E,V^*)$ is generated and $H^0(E^*)=0$.
Hence, by Proposition \ref{p3}, $E$ is semistable and, by
Proposition \ref{p21}, $(E,V^*)$ is $\alpha$-stable for all
$\alpha>0$. So $U^s(n,d,n+1)\ne\emptyset$.

(2) If $d<g+n$ or if $d=g+n$ and $n\not|g$, the bundles $E$
constructed in (1) are stable by Proposition \ref{p3}; hence
$U(n,d,n+1)\ne\emptyset$. If $d=g+n$ and $n|g$, and $g=an$ and
$d=(a+1)n$, in \cite{bu} Butler proves that $E$ is stable unless
$L$ has the form $L\cong L'(Z)$ where $Z$ is an effective divisor
of degree $a+1$ and $L'$ a line bundle with $h^0(L')=n$.

The Brill-Noether number $\beta(1,(a+1)(n-1),n)=0$, hence there
are finitely many choices for $L'$. The dimension of the family
formed of the $L'(Z)$ has dimension $a+1$. Since $a+1<an=g$,
except for $g=n=2$, we can find $L$ lying outside this family. If
$V\subset H^0(L)$ has dimension $n+1$ and generates $L$ then the
kernel of the evaluation map
$$0\rightarrow E^*\rightarrow V\otimes \mathcal{O}\rightarrow
L\rightarrow 0$$ together with the space $V^*$ defines the
generated coherent system $(E,V^*)$ with $E$ stable. By
Proposition \ref{p21}, $(E,V^*)$ is $\alpha$-stable  for all
$\alpha>0$, so $U(n,d,n+1)\ne\emptyset$.
\end{proof}

\begin{lemma}\label{newlemma}
Suppose that $(E,V)\in G(\alpha:n,d,n+1)$ is generically
generated. Then $G(\alpha:n,d,n+1)$ is smooth of dimension $\beta$
at $(E,V)$.
\end{lemma}
\begin{proof} Let $L$ denote the dual of the kernel of the
evaluation map $V\otimes \mathcal{O}\rightarrow E$. The kernel of
the Petri map
\begin{equation}\label{pre1}
V\otimes H^0(E^*\otimes K)\rightarrow H^0(E\otimes E^*\otimes K)
\end{equation}
is $H^0(L ^*\otimes E^*\otimes K).$ Since $E$ is generically
generated from the dual of the exact sequence (\ref{gge}) we have
\begin{equation}\label{e33}
0\rightarrow E^*\otimes L^*\otimes K \rightarrow I^*_E\otimes
L^*\otimes K \rightarrow \tau \rightarrow 0.
\end{equation}

However, since $E$ is generically generated, $I_E$ is generated
and we have the following exact sequence
\begin{equation}\label{ggg1}
0\rightarrow I^*_E\otimes L^*\otimes K\rightarrow V^* \otimes
L^*\otimes K \rightarrow K\rightarrow 0.
\end{equation}

The injectivity of the Petri map for line bundles gives
$H^0(I^*_E\otimes L^*\otimes K)=0$ and from (\ref{e33}),
 $H^0(E^*\otimes L^*\otimes K)=0.$
Therefore, $G(n,d,n+1)$ is smooth of dimension $\beta \geq 0.$
\end{proof}

It is well known that for $n=1,$ the concept of stability is
independent of $\alpha$ and $G(1,d,k):= G(\alpha :1,d,k)=
\mathcal{G}^{k-1}_d$, where $\mathcal{G}^{k-1}_d$ parameterizes
linear series of degree $d$ and dimension
 $k$ (\cite[Chp. 5]{arb}).

 Therefore we have Theorem $2$

\begin{theorem}\label{t21} Let $X$ be general and $ \beta =
\beta(n,d,n+1)\leq g $. Then
\begin{enumerate}
\item $G(\alpha:n,d,n+1)\not= \emptyset$ if and only if $\beta
\geq 0$; \item if $\beta \geq 0$ then $G(n,d,n+1):= G(\alpha
:n,d,n+1) = G(\alpha ':n,d,n+1)$ for all $\alpha ,\alpha ' >0$ and
$\alpha _L =0;$  \item if $\beta> 0$ then $G(n,d,n+1)$ is smooth
and irreducible of dimension $ \beta $ and the generic element is
generated; \item  $U^s(n,d,n+1)= G(n,d,n+1)$ and is birationally
equivalent to $\mathcal{G}^{n}_d$;\item if $\beta =0$
$G(n,d,n+1)\cong \mathcal{G}^{n}_d$ and the number of points of
$G(n,d,n+1)$ is
$$g!\prod_{i=0}^n\frac{i!}{(g-d+n+i)!}.$$
\end{enumerate}
\end{theorem}

\begin{proof} (1) follows from Theorem \ref{t1}(1) and Proposition \ref{p22}.

(2) follows from Propositions \ref{newcor1} and \ref{p7}.

For (3), smoothness follows from Proposition \ref{p7} and  Lemma
\ref{newlemma}. Assume $\beta >0.$ The set of coherent systems
$(E,V)\in G(\alpha :n,d,n+1)$ that are generated is parameterized
by an irreducible variety and has dimension $\beta $ (it is in
correspondence with an open dense set in $B(1,d,n+1), $ which is
irreducible). As in \cite[Theorem 5.11]{bomn}, the irreducibility
of $G(n,d,n+1)$ follows from the fact that the variety that
parameterizes strictly generically generated coherent systems
 has dimension
$< \beta$, so it can not define a new component (see Remark
$2.2$). Hence, $G(n,d,n+1)$ is irreducible.

(4) follows from Proposition \ref{p7} and $(3)$.

For (5), if $\beta=0$, every $(E,V)\in G(n,d,n+1)$ is generated,
hence $G(n,d,n+1)\cong \mathcal{G}^n_d$ which has cardinality
$$g!\prod_{i=0}^n\frac{i!}{(g-d+n+i)!}$$ (see \cite[Chapter V, Theorem 4.4]{arb}).
\end{proof}

\begin{remark}\begin{em} In our case, except when $g=n=2$ and $d=4,
$ Conditions 11.3
in \cite{bomn} are satisfied for $(n,d,n+1),$ i.e. $\beta
(n,d,n+1)\leq n^2(g-1)$, $G_0(n,d,n+1)$ is irreducible and
$B(n,d,n+1)\not= \emptyset$ and hence the results in \cite{bomn}
that assume Conditions 11.3 hold.
\end{em}\end{remark}

\begin{corollary}\label{c21} If $X$ is a general curve and
$0\leq \beta(n,d,n+1)\leq g$ the Brill-Noether locus $B(n,d,n+1)$
is non-empty and irreducible except possibly when $g=n=2$ and
$d=4$. Moreover $G(\alpha:n,d,n+1)$ is a desingularisation of (the
closure of) $B(n,d,n+1).$ The natural map $\phi :
G(\alpha:n,d,n+1) \rightarrow \widetilde{B}(n,d,n+1)$ is an
isomorphism on $B(n,d,n+1)-B(n,d,n+2).$
\end{corollary}

Note that the degree of the bundle $E$ in such coherent systems
satisfies the following inequalities
\begin{equation}
g+n-\frac{g}{n+1} \leq d \leq g+n.
\end{equation}

\begin{proposition}\label{p23}
If $X$ is general and $g\ge n^2-1$, then, for any degree $d\ge
g+n-\frac{g}{n+1}$, $U(n,d,n+1)\ne\emptyset$.
\end{proposition}

\begin{proof} From Proposition \ref{p22} there exist generated
coherent systems $(E,V)$
with $E$ stable for $g+n-\frac{g}{n+1}\leq d \leq g+n.$ Moreover
they are $\alpha $-stable for all $\alpha
>0.$
Given such a coherent system $(E,V)$ and an effective line bundle
$L$ choose a section $s$ of $L$ and define the coherent system
$(E',V')$ as $E':= E\otimes L$ and $V'$ the image of $V$ in
$H^0(E\otimes L)$ under the canonical inclusion $H^0(E)
\hookrightarrow H^0(E\otimes L)$ induced by $s$. It is well known
that $E$ is stable if and only if $E'$ is stable. Moreover, (see
\cite[Lemma 1.5]{rag}) $(E,V)$ is $\alpha $-stable if and only if
$(E', V')$ is $\alpha $-stable.

Therefore, if $g\geq n^2-1$, the length of the interval
$[\frac{g}{n+1}]$ is greater than or equal to $n-1$, so after
tensoring by an effective line bundle, we can obtain all the
values of $d\geq g+n-\frac{g}{n+1}.$
\end{proof}

Moreover, from Theorem \ref{t21}, Proposition \ref{p23}, Lemma
\ref{newlemma} and \cite[Theorem 5.11]{bomn} we have

\begin{theorem}\label{t23} If $X$ is general and $g\geq n^2-1$ then for
any degree $d\geq g+n-\frac{g}{n+1}$
\begin{enumerate}
\item $G(\alpha:n,d,n+1)\not= \emptyset$ for all $\alpha > 0$;
\item $U^s(n,d,n+1)\not= \emptyset $ and is smooth and
irreducible; \item $U(n,d,n+1)\ne\emptyset$ and is smooth and
irreducible.
\end{enumerate}
\end{theorem}

Besides the known relation between coherent systems and
Brill-Noether theory, our results on $G(n,d,n+1)$ can be related
with other problems. Given a generated linear system $(L,V)$ we
have the natural map $$\phi _V :X\rightarrow \mathbb{P}(V^*).$$ In
particular, if $L$ has degree $d$ and dim$V=n+1$ we have;

\begin{theorem}\label{t3} Let $X$ be general, $0
\leq \beta (n,d,n+1)$ and $T\mathbb{P}$ the tangent bundle of
$\mathbb{P}(V^*).$ If $\beta<g$ or $\beta=g$ and $n\not|g$, then
$\phi_V^*(T\mathbb{P})$ is stable. If either $g\ge n^2-1$ or
$\beta=g$, $n|g$ and $g$ and $n$ are not both equal to $2$, then
there exist linear systems $(L,V)$ such that
$\phi_V^*(T\mathbb{P})$ is stable.
\end{theorem}

\begin{proof} Under the hypothesis of the theorem, there exist generated
linear systems $(L,V)$. Denote by $E$ the dual of the kernel of
the evaluation map. Consider the dual Euler sequence
\begin{equation}\label{euler}
0\rightarrow \Omega ^1 _{\mathbb{P}}(1)\rightarrow V\otimes
\mathcal{O}_\mathbb{P}\rightarrow
\mathcal{O}_\mathbb{P}(1)\rightarrow 0
\end{equation}
where $\Omega ^1 _{\mathbb{P}}= T\mathbb{P}^*.$

From the pull-back of (\ref{euler}) we have that $E\otimes
L\cong\phi^*_V(T\mathbb{P})$ (see \cite{gh}). Recall that if $E$
is stable, $E\otimes L$ is stable.

If $\beta<g$ or $\beta=g$ and $n\not|g$, all such $E$ are stable
by the proof of Proposition \ref{p22}. If $\beta=g$, $n|g$ and $g$
and $n$ are not both equal to $2$, some such $E$ are stable, again
by the proof of Proposition \ref{p22}. Finally, if $g\ge n^2-1$,
$U(n,d,n+1)$ is non-empty and irreducible by Theorem \ref{t23} and
its generic element $(E,V^*)$ is generated by the proof of
\cite[Theorem 5.11]{bomn}. Now define $(L,V)$ by dualising the
evaluation sequence of $(E,V^*)$.
\end{proof}

\section{Dual Span}

For a generated coherent system $(E,V)$ of type $(n,d,k)$ with
$H^0(E^*)=0$ denote by $D(E)$ the dual of the kernel of the
evaluation map. That is, we have the following exact sequences
\begin{equation}\label{e31}
0\rightarrow D(E)^*\rightarrow V\otimes \mathcal{O} \rightarrow E
\rightarrow 0
\end{equation}
\begin{equation}
0\rightarrow E^*\rightarrow V^*\otimes \mathcal{O} \rightarrow
D(E) \rightarrow 0
\end{equation}

In \cite[5.4]{bomn}, the coherent system $(D(E),V^*)$ is called
the dual span of $(E,V)$. Note that $(D(E),V^*)$ is a generated
coherent system of type $(k-n,d,k).$ We will define the dual span
for generically generated coherent systems.

 Let $(E,V)$ be a generically generated
coherent system of type $(n,d,k)$ with $H^0(I_E^*)=0$.
  From \cite[Proposition 4.4]{bomn} we have the exact sequence
\begin{equation}\label{1}
  0 \rightarrow N \rightarrow V\otimes \mathcal{O}\rightarrow E
  \rightarrow \tau \rightarrow 0
\end{equation}
with $H^0(N)=0$ and $\tau $ a torsion sheaf of length $\ell$. From
(\ref{1}) we have the exact sequences

\begin{equation}\label{2}
  0 \rightarrow N \rightarrow V\otimes \mathcal{O}\rightarrow I_E
  \rightarrow 0
\end{equation}
 and
 \begin{equation}\label{3}
  0 \rightarrow I_E \rightarrow E
  \rightarrow \tau \rightarrow 0 .
\end{equation}

\begin{lemma}\label{lemma0} $N=D(I_E)^*$.
\end{lemma}
\begin{proof} The coherent system $(I_E,V)$ is generated. From (\ref{2})
 $N=D(I_E)^*$.
\end{proof}

\begin{remark}\begin{em}\label{remark00}
 If $(E,V)$ is generically generated and $H^0(I_E^*)=0$, from (\ref{1}) and Lemma
\ref{lemma0} we have the sequence
\begin{equation}\label{4t}
  0 \rightarrow D(I_E)^* \rightarrow V\otimes \mathcal{O}\rightarrow E
  \rightarrow \tau \rightarrow 0.
\end{equation}
and
\begin{equation}
0\rightarrow D(I_E)^*\rightarrow V\otimes \mathcal{O} \rightarrow
I_E \rightarrow 0
\end{equation}
Moreover,  $(D(I_E),V^*)$ is the dual span of $(I_E,V)$.
\end{em}\end{remark}

\bigskip

Let
\begin{equation}\label{444}
  0 \rightarrow D(I_E) \rightarrow D(E)_\ell\rightarrow
  \tau ' \rightarrow 0
\end{equation}
be an elementary transformation of $D(I_E)$ with $\tau '$ a
torsion sheaf of length $\ell.$ The subspace $V^*\subset
H^0(D(I_E))$ defines a subspace $V'$ in $H^0(D(E)_\ell)$ which we
identify with $V^*$.

\begin{definition}\label{def1} Let $(E,V)$ be a generically
generated coherent system of type $(n,d,k)$ with $H^0(I_E^*)=0$. A
dual span of $(E,V)$, denoted by $D(E,V)$, is an elementary
transformation $(D(E)_\ell ,V^*)$ of $(D(I_E),V^*)$ of length
$\ell $ where $\ell =d_E - d_{{I_E}}.$
\end{definition}

\begin{remark}\begin{em}\label{remark1}
\begin{enumerate} \item If $(E,V)$ is strictly generically
generated then the family of dual spans associated to $(E,V)$ has
dimension at most $\ell n -1$; \item If $(E,V)$ is generated there
is a unique dual span given by $(D(E),V^*)$; \item If $(E,V)$ is a
generically generated coherent system of type $(n,d,k)$,
$(D(I_E),V^*)$ is a generated coherent system of type
$(k-n,d-\ell,k);$ \item $D(E,V)$ is a coherent system of type
$(k-n,d,k)$;\item  the image of the evaluation map $V^*\otimes
\mathcal{O} \rightarrow D(E)_\ell $ is $D(I_E).$
\end{enumerate}
\end{em}\end{remark}

\begin{proposition}\label{propos} Let $(E,V)$ be a coherent systems of type
$(n,d,k)$. If $(E,V)$ is generically generated with $H^0(I^*_E)=0$
then a dual span $D(E,V)=(D(E)_\ell,V^*)$ is generically
generated. Moreover, $H^0(I^*_{{D(E)_\ell}})=0.$
\end{proposition}

\begin{proof} The Proposition follows from the definition of a dual
span, since $(D(I_E),V^*)$ is generated and
$I_{{D(E)_\ell}}=D(I_E)$.
\end{proof}

\begin{remark}\begin{em} Note from the definition of a dual span
that $(E,V)$ is a dual span of $D(E,V)=(D(E)_\ell,V^*)$.
\end{em}\end{remark}

\begin{theorem}\label{t31} Let $X$ be a general
 curve of genus $g$ and $d< g+n_1\leq g+n_2$ then for all $\alpha
 >0$, $G(\alpha:n_1,d, n_1+n_2)\not= \emptyset $ if and only if
$G(\alpha:n_2,d, n_1+n_2)\not= \emptyset $.
\end{theorem}

\begin{proof} Let $(E,V)\in G(\alpha:n_i,d,n_1+n_2)$ for $i=1,2$. From
Proposition \ref{p7}, $(E,V)$ is generically generated and
$H^0(I^*_E)=0$. From Proposition \ref{propos} a dual span
$D(E,V)=(D(E)_\ell,V^*)$ is generically generated with
$H^0(I^*_{{D(E)_\ell}})=0$ and from Proposition \ref{p6} it is
$\alpha $-stable for all $\alpha >0$.
\end{proof}

For any $(n,d,k)$ define $G_g$ as
$$G_g(n,d,k):=\{(E,V): (E,V) \mbox{ is of type } (n,d,k) \ \mbox{ and it is
generated with $H^0(E^*)=0$}\}$$

\begin{corollary}\label{teo1} If $d<g+n_1\leq g+n_2$, then
 $G_g(n_i,d,n_1+n_2)\subset U(n_i,d,n_1+n_2)$ for
$i=1,2$. Moreover, for $i=1,2$, $G_g(n_i,d,n_1+n_2)$ is open and
$G_g(n_1,d,n_1+n_2)\cong G_g(n_2,d,n_1+n_2).$
\end{corollary}

\begin{proof} If $(E,V)\in G_g(n_i,d,n_1+n_2)$, from Proposition
\ref{p6} is $\alpha $-stable for all $\alpha >0$ and from
Proposition \ref{p3} $E$ is stable. The dual span correspondence
for generated coherent systems gives the isomorphism.
\end{proof}

 To prove Theorem $7$ we give
four Lemmas that we will use

\begin{lemma}\label{lem1} Let $(E,V)$ be a generated coherent
system. The Petri map of $(E,V)$ is injective if and only if the
Petri map of $D(E,V)$ is injective.
\end{lemma}

\begin{proof} Since $(E,V)$ is generated $D(E,V)=(D(E),V^*)$.
We have the following exact sequences
\begin{equation}
0\rightarrow D(E)^*\rightarrow V\otimes \mathcal{O} \rightarrow
E\rightarrow 0,
\end{equation}
\begin{equation}
0\rightarrow E^*\rightarrow V^*\otimes \mathcal{O} \rightarrow
D(E)\rightarrow 0.
\end{equation}
The Lemma follows from the cohomology sequences
\begin{equation}
0\rightarrow H^0(D(E)^*\otimes E^*\otimes K) \rightarrow V\otimes
H^0(E^*\otimes K) \stackrel{\psi}{\rightarrow} H^0(E\otimes
E^*\otimes K)\cdots
 \end{equation}
\begin{equation}
0\rightarrow H^0(E^*\otimes D(E)^*\otimes K)\rightarrow V\otimes
H^0(D(E)^*\otimes K) \stackrel{\phi}{\rightarrow} H^0(D(E)
 \otimes D(E)^*\otimes K)\cdots
 \end{equation}
since $\phi $ is injective if and only if $\psi $ is injective.
\end{proof}

\begin{lemma}\label{lem2} Let $(E,V)$ be strictly generically
generated. If the Petri map of $(I_E,V)$ is injective, the Petri
map of $(E,V)$ is injective.
\end{lemma}

\begin{proof} The Lemma follows from the cohomology sequences
\begin{equation}
0\rightarrow H^0(D(I_E)^*\otimes I_E^*\otimes K) \rightarrow
V\otimes H^0(I_E^*\otimes K) \stackrel{\psi}{\rightarrow}
H^0(I_E\otimes I_E^*\otimes K)\cdots
 \end{equation}
 \begin{equation}
0\rightarrow H^0(D(I_E)^*\otimes E^*\otimes K) \rightarrow
V\otimes H^0(E^*\otimes K) \stackrel{\psi}{\rightarrow}
H^0(E\otimes E^*\otimes K)\cdots
 \end{equation}

and the cohomology of the exact sequence

\begin{equation}\label{ee31}
  0 \rightarrow E^*\otimes D(I_E)^*\otimes K \rightarrow
  I^*_E\otimes D(I_E)^*\otimes K
  \rightarrow \tau \rightarrow 0 .
\end{equation}

\end{proof}

Let $(E,V)$ be a generically generated coherent system. From
Proposition \ref{propos} a dual span $D(E,V)=(D(E)_\ell,V^*)$ is
generically generated. Hence, from Remark \ref{remark00}, we have
the sequence
\begin{equation}\label{41}
  0 \rightarrow I_E^* \rightarrow V^*\otimes
  \mathcal{O}\rightarrow D(E)_\ell
  \rightarrow \tau \rightarrow 0
\end{equation}

\begin{lemma}\label{lemmaf} The Petri map of $(I_E,V)$ is
injective if and only if the Petri map of $(I_{{D(E)_\ell}},V^*)$
is injective.
\end{lemma}

\begin{proof} The Lemma follows at once from Lemma \ref{lem1}
since $I_{{D(E)_\ell}}= D(I_E).$
\end{proof}

\begin{lemma}\label{lem3} If the Petri map of $(I_E,V)$ is
injective, the Petri map of a dual span $D(E,V)=(D(E)_\ell,V^*)$
is injective.
\end{lemma}

\begin{proof} From (\ref{41}), the kernel
of the Petri map of $(D(E)_\ell,V^*)$ is $H^0(I_E^*\otimes
D(E)^*_\ell \otimes K).$

From the exact sequence (\ref{444}) we obtain the following exact
sequence
\begin{equation}\label{4}
0\rightarrow D(E)^*_\ell \otimes I^*_E\otimes K \rightarrow
D(I_E)^*\otimes I^*_E \otimes K
  \rightarrow \tau \rightarrow 0.
\end{equation}

The kernel of the Petri map for $(I_E,V)$ is $H^0(D(I_E)^*\otimes
I_E^*\otimes K).$ Hence, if $H^0(D(I_E)^*\otimes I^*_E\otimes
K)=0$, $H^0(I_E^*\otimes D(E)^*_\ell \otimes K)=0.$
\end{proof}

We now have Theorem $7$

\begin{theorem}\label{theorem1} Let $(E,V) \in
G(\alpha:n_1,d,n_1+n_2)$. If either of the Petri maps of $(I_E,V)$
or $(I_{{D(E)_\ell}},V^*)$ is injective then,
\begin{enumerate} \item $G(\alpha:n_1,d,n_1+n_2)$ is smooth of dimension
$\beta(n_1,d,n_1+n_2)$ in a neighbourhood of $(E,V)$; \item
$G(\alpha:n_2,d,n_1+n_2)$ is smooth of dimension
$\beta(n_2,d,n_1+n_2)$ in a neighbourhood of the dual span
$D(E,V).$
\end{enumerate}
\end{theorem}

\begin{proof} If the Petri map of $(I_E,V)$ is injective, from Lemmas \ref{lem1},
  \ref{lem2} and \ref{lem3}, the Petri maps of $(E,V)$ and $D(E,V)$ are injective.
 From Proposition \ref{e} $G(\alpha:n_i,d,n_1+n_2)$
 is smooth of dimension $\beta(n_i,d,n_1+n_2)$ in a neighbourhood of
 $(E,V)$ and of $D(E,V)$, respectively.

If the Petri map of $(I_{{D_E}},V)$ is injective, again from
Lemmas \ref{lem1} and \ref{lem2}, the Petri map of $D(E,V)$ is
injective. From Lemma \ref{lemmaf}, the Petri map of $(I_E,V)$ is
injective and, as above, the Petri map of $(E,V)$ is injective.
Hence, $G(\alpha:n_i,d,n_1+n_2)$
 is smooth of dimension $\beta(n_i,d,n_1+n_2)$ in a neighbourhood of
 $D(E,V)$ and of $(E,V)$, respectively.

\end{proof}

\begin{remark}\begin{em} Theorems \ref{t31} and \ref{theorem1}
 apply for any
$\alpha >0.$ Since $d<g+n_1$ the bundles in $G(\alpha:
n_1,d,n_1+n_2)$ are stable (see Proposition \ref{p7}). Hence, we
have similar results for the Brill-Noether loci
${B}(n_1,d,n_1+n_2)$ and ${B}(n_2,d,n_1+n_2)$.
\end{em}
\end{remark}

\section{ Rank 2 and Genus 2}

In this section we will consider the case $n=2$ and then $g=2.$

From Proposition \ref{p23} we have that for a general curve and
$g\geq 3$, $G(\alpha ;2,d,3)\not= \emptyset$
 for all $\alpha >0$ and $U(2,d,3)\not= \emptyset$
 for $d\geq \frac{2g}{3} +2$. For
 $k>4$, we have the following Theorem

\begin{theorem}\label{teo7} Let $X$ be general, $s\geq 3$ and
$d< s +2g -\frac{4g}{s+2}$. If $G_0(2,d,2+s)$ is non-empty then
$G(\alpha :2,d,2+s)$ is non-empty for all $\alpha >0.$ Moreover,
$U^s(2,d,2+s)\not= \emptyset.$
\end{theorem}

\begin{proof} Let $(E,V)\in G_0(2,d,2+s)$. Hence, $E$ is semistable.

Let $r_s :=\lceil\frac{2+s}{2}\rceil $ and
 $(F,W)$ a coherent subsystem of $(E,V)$ with $n_F=1$. If $\dim W \geq r_s$ the
Brill-Noether number $\beta (1,d_F,r_s)\geq 0$. That is, $ d_F\geq
r_s+g-1 -\frac{g}{r_s}$. But then
$$ d_F\geq r_s+g-1 -\frac{g}{r_s}> \frac{d}{2},$$
which is a contradiction since $E$ is semistable. Therefore, for
any coherent subsystem $(F,W)$ $\dim W < \frac{2+s}{2}$ and from
Remark  \ref{r0} $(3)$ $(E,V)$ is $\alpha $-stable for all $\alpha
>0.$ Therefore, $G(\alpha :2,d,2+s)\not = \emptyset$
for all $\alpha >0$ and  $U^s(2,d,2+s)\not= \emptyset.$
\end{proof}

Let $X$ be any curve. From Proposition \ref{p21} we have that any
generated coherent system $(E,V)$ of type $(n,d,n+1)$ with $E$
stable is $\alpha $-stable for all $\alpha >0$. For $n=2$ we have
(see \cite[Theorem 9.2]{bomn} for general curve).

\begin{proposition}\label{propn2} Let $X$ be any curve. If
$G_0(2,d,4)\not = \emptyset$ and there exists a generated coherent
system $(E,V)\in G_0(2,d,4)$, then ${G}(\alpha:2,d,4)\not=
\emptyset$ for all $\alpha >0$ and $U^s(2,d,4)\not= \emptyset$.
Moreover, if $E$ is stable $U(2,d,4)\not= \emptyset .$
\end{proposition}

\begin{proof}
Let $(F,W)$ be a coherent subsystem of $(E,V)$ with $n_F=1$. From
Proposition \ref{p21}, $\dim W \leq 2$. If $\dim W=2$, since
$(E,V)\in G_0(2,d,4)$, $d_F<\mu(E)$. From Remark \ref{r0} $(E,V)$
is $\alpha$-stable for all $\alpha >0.$
\end{proof}

\begin{corollary}\label{t1n2} For any curve $X$ and $d\geq 4g-2$,
$G(\alpha:2,d,4)\not= \emptyset$ for all $\alpha >0$. Moreover,
$U(2,d,4)\not= \emptyset$ and for $d\geq 4(g-1)$, $U(2,d,2+s)=
\emptyset$ if $s>d-2g.$
\end{corollary}

\begin{proof} Since any stable bundle of degree $d\geq 2(2g-1)$ is
generated, the first part follows from Proposition \ref{propn2}.
The last part follows from the Riemann-Roch Theorem.
\end{proof}

\begin{remark}\begin{em}\label{r31}
Recall from the Brill-Noether theory for vector bundles of rank
$n\geq 2$ (see \cite{bgn}, \cite{me} and \cite{bfmn}) that if
$0<d< 2n$, there exists a semistable vector bundle $E$ of rank $n$
and degree $d$ with $k$ sections if and only if $n\le d+(n-k)g$.
Hence, if $0<d<2n$ and $k>n+\frac{d-n}g$, then
$U^s(n,d,k)=\emptyset$. Moreover, if $d>n(2g-2)$, then by the
Riemann-Roch theorem every semistable bundle $E$ has
$h^0(E)=d+n(1-g)$; so, if $k>d+n(1-g)$, $U^s(n,d,k)=\emptyset$.
\end{em}
\end{remark}

We shall now consider the case $g=2$. Any curve of genus $g=2$ is
a Petri curve. From Corollary \ref{corfinal}, if $d<n+2,$
$G(\alpha:n,d,k)=\emptyset$ for all $\alpha >0$ and
$k>n$.

From Theorem \ref{t21} we have

\begin{proposition}\label{c41} For $X$  of genus $g=2$ and $d=n+2$, $n\geq
3$;
\begin{enumerate}
\item $G(\alpha:n,d,n+1)\not= \emptyset $ for all $\alpha >0$;
\item $G(\alpha:n,d,n+1)= G(\alpha ':n,d,n+1)$ for all $\alpha ,
\alpha '>0$ and $\alpha _L =0;$
 \item $G(n,d,n+1)$ is smooth and irreducible of dimension $2$;
\item $U(n,d,n+1)= G(n,d,n+1)$; \item if $k\geq n+2,$
$G(\alpha:n,d,k)= \emptyset$ for all $\alpha >0.$
\end{enumerate}
\end{proposition}

\begin{proof} Parts 1), 2), 3) and 4) follow from Theorem \ref{t21}.
Part 5) follows from Remark \ref{r31} and Proposition $3.7$, since
for the existence of a semistable bundle with at least $k$
sections we need $k-n\leq \frac{d-n}{2}.$
\end{proof}

From Remark \ref{r31} and Proposition \ref{p1dg} we have
\begin{enumerate}
\item if $n+2<d<2n$ and $k>\frac{d+n}2$, $U^s(n,d,k)=\emptyset$;
\item if $d>2n$ and $k>d-n$, $U^s(n,d,k)=\emptyset$; \item if
$d\geq 3n$, $G(\alpha:n,d,n+1)\not= \emptyset$ for all $\alpha
>0$. Moreover, $U(n,d,n+1)\not= \emptyset .$
\end{enumerate}

In particular for $n=2$, from Propositions \ref{p21}, and
\ref{t1n2}, Corollary \ref{corfinal} and the Riemann-Roch Theorem
we have

\begin{enumerate}
\item If $d<4$ and $k\geq 3$, $G(\alpha :2,d,k)=\emptyset$ for all
$\alpha
>0$;
 \item if $d= 5$ and $k> 3$, $U(2,d,k)=\emptyset $ and $G_0(2,5,k)=\emptyset$; \item if
$d\geq 6$ and $k=3,4$,  $G(\alpha :2,d,k)\not =\emptyset$ for all
$\alpha >0$. Moreover, $U(2,d,k)\not=\emptyset ;$ \item if $d \geq
6$ and $k> d-2$, $U(2,d,k)=\emptyset $ and $G_0(2,d,k)=\emptyset
.$
\end{enumerate}

 For $d=4$ we need the following Lemmas

\begin{lemma}\label{corn2} \begin{enumerate}\item $B(2,4,k)=\emptyset
$ for $k\geq 3$; \item $\widetilde{B}(2,4,k)=\emptyset$ for $k\geq
5$; \item $\widetilde{B}(2,4,3)\not= \emptyset$; \item
$\widetilde{B}(2,4,4)=\{K\oplus K\}$.
\end{enumerate}
\end{lemma}

\begin{proof}
 Let $E$ be a semistable vector bundle of rank $2$ and
degree $d=4=2(2g-2)$. From the Riemann-Rock theorem,
$h^0(E)=2+h^1(E)$. If $h^1(E)= h^0(E^*\otimes K)\geq 1 $, then $E$
is an extension

\begin{equation}\label{eq12}
\xi :0\rightarrow L\rightarrow E \rightarrow K\rightarrow 0
\end{equation}
of $K$ by $L$, where $L$ is a line bundle of degree $2$. Thus, $E$
can not be stable. That is, $B(2,4,k)= \emptyset $ for $k\geq 3$.

Since $h^1(L)\leq 1$ and $h^1(K)=1$, from the cohomology sequence
of (\ref{eq12}), $h^1(E)\leq 2.$ Hence, $\widetilde{B}(2,4,k)=
\emptyset$ for $k \geq 5.$

If $L\not\cong K$, $H^1(L)=0$, $H^0(L)\cong \mathbb{C}$ and
$h^1(K^*\otimes L)=1$. Hence, there exist non-trivial extensions
(\ref{eq12}), and $h^0(E)=3$. That is, $\widetilde{B}(2,4,3)\not=
\emptyset .$

Let $L\cong K$. If $\xi $ is non-trivial, from the cohomology
sequence of $$0\rightarrow \mathcal{O}\rightarrow E^*\otimes K
\rightarrow \mathcal{O}\rightarrow 0,$$ $H^0(E^*\otimes K)\cong
H^0(\mathcal{O})$. Hence, $h^0(E)=3$.

Therefore, $\widetilde{B}(2,4,4)=\{K\oplus K\}.$
\end{proof}

Note that if $(L,W)$ is a coherent system of type $(1,2,2)$ then
$(L,W)=(K,H^0(K)).$

\begin{lemma}\label{ppp} If $(K,H^0(K))$ is a coherent
subsystem of a coherent system $(E,V)$ of type $(2,4,3)$ then
$(E,V)$ is not $\alpha $-semistable for any $\alpha >0$.
\end{lemma}

\begin{proof} For any $\alpha >0$, $\mu _{\alpha
}(K,H^0(K))=2+2\alpha > 2+\alpha \frac{3}{2}=\mu _{\alpha }(E,V).$
\end{proof}

\begin{corollary}\label{corff} The coherent systems

\begin{enumerate} \item $(L\oplus K, H^0(L)\oplus H^0(K));$ \item
$(E,H^0(E))$ with $E$ a non-trivial extension of $K$ by $K;$
\end{enumerate}
are not $\alpha $-semistable for any $\alpha >0$.
\end{corollary}

\begin{lemma}\label{propn2f} Let $(E,V)$ be a coherent system
of type $(2,4,3)$. If $E$ is a non-trivial extension $\xi $ of $K$
by $L$, with $L\not\cong K$, $(E,V)$ is generated. Moreover,
$(E,V)$ is $\alpha$-stable for all $\alpha >0$.
\end{lemma}

\begin{proof} If $n_{I_E}=1$ then $I_E=K$, which is a contradiction
since $\xi \not= 0$. If $n_{I_E}=2$ and $d_{I_E}<4$, from Lemma
\ref{re1} we get a contradiction. Therefore, $(E,V)$ is generated.
From Proposition \ref{p21}, $(E,V)$ is $\alpha $-stable for all
$\alpha >0$.
\end{proof}

\begin{proposition}\label{final}
If $(E,V)\in G_L(2,4,3)$ then $E$ is semistable and $(E,V)$ is
$\alpha $-stable for all $\alpha >0$.
\end{proposition}

\begin{proof} The Proposition follows at once from Propositions \ref{newcor1} and \ref{p7}.
\end{proof}

\begin{theorem}\label{proposition4} \begin{enumerate}\item $U(2,4,k)=\emptyset
$ for $k\geq 3$; \item $G_0(2,4,k)=\emptyset$ for $k\geq 5$; \item
$G(\alpha :2,4,3)\not= \emptyset $ for all $\alpha >0$; \item
$U^s(2,4,3)\cong G_L(2,4,3)\cong Pic^4(X)$.
\end{enumerate}
\end{theorem}

\begin{proof}$(1)$, $(2)$ and $(3)$ follow from Lemma \ref{corn2}.
$U^s(2,4,3)\cong G_L(2,4,3)$,  follow from Lemma \ref{propn2f} and
Proposition \ref{final}.

To prove $G_L(2,4,3)\cong Pic^4(X)$ suppose $(E,V)\in
G_L(2,4,3)\cong U^s(2,4,3)$, so $E$ is semistable, generically
generated and $h^0(I_E^*)=0$. If $E$ is not generated, then $\deg
I_E\leq 3$. Moreover $I_E$ must be stable, for otherwise it has a
quotient line bundle $Q$ of degree $\le1$, hence with $h^0(Q) \leq
1$. The corresponding subbundle $L$ has $\dim(V\cap H^0(L))\ge2$,
contradicting the $\alpha$-stability of $(E,V)$. However
$U(2,3,3)=\emptyset$, so $I_E$ cannot exist. Thus $E$ is generated
and it follows that $E$ arises from an extension
$$0\rightarrow L^*\rightarrow V\otimes \mathcal{O}\rightarrow E\rightarrow 0$$
or dually
\begin{equation}\label{new}
0\rightarrow E^*\rightarrow V^*\otimes \mathcal{O}\rightarrow
L\rightarrow0,
\end{equation}
where $L$ is a line bundle of degree $4$.

 Conversely, any line
bundle $L$ of degree $4$ is generated and $h^0(L)=3$ by the
Riemann-Roch Theorem. So there is a unique extension (\ref{new})
for each $L$. Certainly then $E$ is generated with $h^0(E^*)=0$,
so $(E,V)\in G_L(2,4,3)$.
\end{proof}

Moreover,

\begin{theorem}\label{teoff}
\begin{enumerate}
\item $\widetilde{G}_0(2,4,4)=\{(K\oplus K , H^0(K\oplus K))\}$;
\item $G_0(2,4,4)= \emptyset$;  \item
$\widetilde{G}(\alpha:2,4,4)\not= \emptyset$ for all $\alpha >0$.
\end{enumerate}
\end{theorem}

\begin{proof} Parts $(1)$ and $(2)$ follow from Lemma \ref{corn2}.
Since $(K\oplus K, H^0(K\oplus K))\cong
(K,H^0(K))\oplus(K,H^0(K)),$ it is $\alpha$-semistable for all
$\alpha >0$, so
 $\widetilde{G}(\alpha:2,4,4)\not= \emptyset$ for all $\alpha >0$.

\end{proof}

For $d=5$ and $k=3$ we have

\begin{theorem}\label{prop5}
\begin{enumerate} \item $G_0(2,5,3)\not= \emptyset$;
 \item $U(2,5,3)\not =\emptyset $; \item $U(2,5,3)\not=
 G_0(2,5,3).$
\end{enumerate}
\end{theorem}

\begin{proof} Let $E$ be a non-trivial
extension
\begin{equation}\label{seq6}
\phi : 0\rightarrow L\rightarrow E\rightarrow M\rightarrow 0
\end{equation}
of $M$ by $L$, where $L $ is a line bundle of degree $2$ and $M$ a
general line bundle of degree $3$ with $h^0(M)=2.$ Note
 that, $h^1(M^*\otimes  L)=2.$

 It is well known that $E$ is stable and from the cohomology
 sequence of $(\ref{seq6})$, $h^0(E)=3.$ Hence, $(E,H^0(E))\in G_0(2,5,3).$

Let $(F,W)$ be any coherent subsystem of $(E,H^0(E))$, with $F$ a
line sub-bundle. Since $E$ is stable, $d_F<\mu (E)=2+\frac{1}{2},$
so $\dim W\leq h^0(F)\leq 2$. If $\dim W=2$, $F\cong K.$

Now, if in the extension $(\ref{seq6})$ $L\not \cong K$ and $M$ is
general and generated, $H^0(K^*\otimes M)=0$ i.e. $K$ can not be a
subbundle of $E$. Hence, for all coherent subsystems $(F,W)$,
$\dim W\leq 1$ and from Remark \ref{r0}, $(E,H^0(E))$ is $\alpha
$-stable for all $\alpha >0$. Therefore, $U(2,5,3)\not=\emptyset
.$

However, if $L\cong K$, for any coherent subsystem $(F,W)$ of
$(E,H^0(E))$, with $F$ a line subbundle, $\mu _{\alpha}(F,W)\leq
\mu _{\alpha } (K,H^0(K))$. Thus, since $(K,H^0(K))$ is a coherent
subsystem of $(E,H^0(E)),$ $$\mu _{\alpha}(K,H^0(K))<\mu _{\alpha
} (E,V)$$ if and only if $\alpha <1$. For $\alpha =1$, $\mu
_{\alpha}(K,H^0(K))=\mu _{\alpha } (E,V).$ Therefore, $(E,H^0(E))
\not \in U(2,5,3).$

\end{proof}

 {\small \noindent \textbf{Acknowledgements.} The author thanks
 Usha Bhosle for useful
conversations and specially Peter Newstead for his suggestions and
comments on previous versions of this work. Thanks are also due to
the referee for valuable comments and suggestions towards
improvement of the paper. She also thanks the International Center
for Theoretical Physics, Trieste, where part the work was carried
out, for their hospitality and acknowledges the support of CONACYT
grant 48263-F. The author is a member of the research group VBAC
(Vector Bundles on Algebraic Curves), which was partially
supported by EAGER (EC FP5 Contract no. HPRN-CT-2000-00099) and by
EDGE (EC FP5 Contract no. HPRN-CT-2000-00101).}

%%%%%%%%%%%%%%%%%%%%%%%%%%%%%%%%%%%%%%%%%%%%%%%%%%%%%%%%%

\end{document}